\documentclass{article}     

\usepackage{amsmath}
\usepackage{amssymb}
\usepackage{amsthm}
\usepackage{bm}

\usepackage{enumerate}

\newtheorem{theo}{Theorem}[section]  
\newtheorem{coro}[theo]{Corollary}
\newtheorem{lemma}[theo] {Lemma} 
\newtheorem{prop}[theo] {Proposition}

\newtheorem{defi}[theo]{Definition}

\begin{document}
	
	\setcounter{page}{1}   
	
	\title{The Finite Model Property of Some Non-normal Modal Logics with the Transitivity Axiom}
	\author{Kirill Kopnev}
	
	\maketitle

	\begin{abstract} 
		
		In 1997 Timothy J. Surendonk proved via algebraic semantics that all modal logics without iterative axioms are canonical and so strongly complete. In this paper we continue the work done by Surendonk in this field. We use neigborhood semantics to show that some iterative logics (with the axiom $\Box p\rightarrow\Box\Box p$) are also strongly complete and have the finite model property. 
		
	\end{abstract}


	\section*{Introduction and overview}

	This paper explores completeness and finite model property (FMP) for some non-normal propositional modal logics with the axiom $(4)~\Box p\rightarrow \Box\Box p$. Non-normal modal logics are logics weaker then well known logic $\textup{K}$, e.g. logics without axioms $\Box(p\rightarrow q)\rightarrow(\Box p\rightarrow\Box q)$ or $\Box\top$. It should be mentioned, that the problem of FMP for logics with the transitivity axiom was posed in \cite{art:Surendonk} by Timothy Surendonk. 
	
	Study of non-normal modal logics is usually motivated by working with some epistemic problems or questions related to multi-agent systems. Since non-normal logics are weaker than $\textup{K}$ and $\textup{K}$ is sound with respect to Kripke frames (cf.\cite{book:BBB}), we have that Kripke semantics is not adequate for non-normal modal logics. Nevertheless, Dana Scott \cite{art:Scott} and Richard Montague \cite{art:Montague} independently introduced another semantics. This semantics is known as neighborhood semantics nowadays. In addition to the fact that it allows us to reason formally about different philosophical problems, neighborhood semantics also might be interesting for mathematicians due to its connection with topology. Briefly this connection can be described in the following way: topological spaces can be viewed as a certain kind of neighborhood structures. In the language of modal logic these are defined by the axioms of the logic $\textup{S4}$. The completeness result for the logic $\textup{S4}$ with respect to topological semantics was proved by John McKinsey and Alfred Tarski in 1944 \cite{art:McKinseyTarski}. 
	
	Moreover, neighborhood semantics can be used to explore very important class of mathematical objects using tools of modal logic. These objects are so-called closure operators. In terms of the modal language this is the logic with the following axioms: $(\textup{T})\;\Box p\rightarrow p, (\textup{M})\;\Box(p\wedge q)\rightarrow(\Box p\wedge\Box q)$ and $(4)$.
	
	Closure operators are used in different fields of mathematics: algebra, logic, topology, etc. In terms of sets, a closure operator on a set $X$ is a function $Cl\colon\rho(X)\rightarrow\rho(X)$ satisfying the following conditions for every $A, B\subseteq X$:

	\begin{itemize}
		\item $A\subseteq Cl(A)$;
		
		\item $A\subseteq B\Rightarrow Cl(A)\subseteq Cl(B)$;
		
		\item $Cl(Cl(A))=Cl(A)$.
		
	\end{itemize}
	
	Note that rewriting axioms of $\textup{S04}$ in terms of $\Diamond$ gives us exactly the definition of a closure operator.
	
	A particular type of closure operators is Tarski's consequence operator. They were introduced by Tarski in \cite{art:Tarski} for investigation of deductive systems. For a set of sentences $A$, the set of all consequences of $A$ is the intersection of all sets containing the set $A$ and are closed under the given rules of inference. This is a special case of a closure operator defined in terms of closed sets. Consult \cite{art:Tarski} to check that Tarski's consequence operator has all three properties listed above. An important property of Tarski's consequence operators is that $\forall X\subseteq A: Cl(X)=\bigcup\limits_{i\in I}Cl(X_{i})$, where $X_{i}\subseteq A$ and $|X_{i}|<\omega$ for each $i\in I$. Closure operators of this kind are finitary. 
	
	In algebra closure operator appears in several contexts. The best known one is linear algebra where the closure operator is understood as the function associating every subset of a given vector space with its linear span. The closure operator is finitary in this case as well. The additional property of this linear span as a closure operator in a vector space is the exchange property: if $x$ is in the closure of the union of $A$ and $\{y\}$ but not in the closure of $A$, then $y$ is in the closure of the union of $A$ and $\{x\}$. A closure operator with this property is called a matroid. 
	
	Our work is, in a sense, a continuation of David Lewis' result from \cite{art:Lewis}, stating that every intensional logic without iterative axioms has the finite model property. It should be mentioned, that Lewis' paper contains an error in the proof of Lemma $11$ while applying Lemma $7$ for the second time (cf. appendix \ref{appendix.A}). Nevertheless, Surendonk \cite{art:Surendonk} proved Lewis' theorem using algebraic semantics. The result by Surendonk is sufficient for establishing the same properties for non-iterative logics in neighborhood semantics. In our paper, we prove the same result for some iterative modal logics using a different technique.

	\section{Language and Logics} %
	
	Let us first recall the definition of modal formulas.
	
	\begin{defi}
		Suppose that $Var=\{p,q,r,...\}$ is a countable set of proposition variables. \textit{Modal formulas} are constructed recursively from set $Var$ using the logical connectives negation $(\neg)$, conjunction $(\wedge)$, and the unary modal operator $\Box$.
	\end{defi}
	
	Additional propositional connectives $(\lor, \rightarrow, \leftrightarrow)$ are defined as usual. We will also introduce special formulas $\top$ and $\bot$ as abbreviations for $\neg(p\wedge\neg p)$ and $(p\wedge\neg p)$ respectively. Modal connective $\Diamond$ is defined in the following way: $\Diamond\varphi:=\neg\Box\neg\varphi$.
	
	Before defining a logic, let us first define substitutions.
	
	\begin{defi}
		A  \textit{substitution} $\sigma$ is a function $\sigma\colon Var\rightarrow Fm$. A substitution $\sigma$ is extended to a formula homomorphism $\overline{\sigma}\colon Fm\rightarrow Fm$ by recursion:
		\\1. $\overline{\sigma}(p)=\sigma(p)$, for all $p\in Var$;
		\\2. $\overline{\sigma}(\neg\psi)=\neg\overline{\sigma}(\psi)$;
		\\3. $\overline{\sigma}(\psi\wedge\gamma)=(\overline{\sigma}(\psi)\wedge\overline{\sigma}(\gamma))$;
		\\4. $\overline{\sigma}(\Box\psi)=\Box\overline{\sigma}(\psi)$;
	\end{defi}
	
	We identify a logic with the set of its theorems.
	
	\begin{defi}
		
		A set of formulas $\textup{L}$ is called a  \textit{modal logic} if $\textup{L}$ is closed under the rules $(MP)\left(\displaystyle{\varphi\rightarrow\psi, \varphi\over \psi}\right), (Sub)\left(\displaystyle{\varphi\over\overline{\sigma}(\varphi)}\right), (RE)\displaystyle\left({{\varphi\leftrightarrow \psi}\over{\Box\varphi\leftrightarrow\Box\psi}}\right)$ and contains all classical tautologies. 
		
	\end{defi} 
	
	As in this paper we consider only modal logics, we omit the word `modal'.
	
	A set of formulas $A$ is an  \textit{axiom set} for a logic $\textup{L}$ iff $\textup{L}$ is the least logic containing $A$. The smallest logic obtained by adding $\Gamma\subseteq Fm$ to $\textup{L}$ is denoted by $\textup{L}+\Gamma$. If $\Gamma=\{\varphi\}$, then we write $\textup{L}+\varphi$ instead of $\textup{L}+\{\varphi\}$.
	
	The smallest logic is denoted by $\textup{E}$.
	
	In our paper, we consider logics with the following extra axioms:
	
	\begin{flushleft}
		$(4)~\Box p\rightarrow\Box\Box p$;
		
		$(\textup{T})\; \Box p\rightarrow p$;
		
		$(\textup{M})\; \Box(p\wedge q)\rightarrow(\Box p\wedge\Box q)$;
		
		$(\textup{C})\; (\Box p\wedge\Box q)\rightarrow\Box(p\wedge q)$.
	\end{flushleft} 
	
	So we will cosider with the following logics:
	
	\begin{flushleft}
		$\textup{E}4=\textup{E}+4$;
		
		$\textup{EMC4}=\textup{E}4+\textup{M}+\textup{C}$;
		
		$\textup{S04}=\textup{E}4+\textup{T}+\textup{M}$.
	\end{flushleft}
	
	Since $\vdash_{\textup{EMC}} \Box(p\rightarrow q)\rightarrow(\Box p\rightarrow\Box q)$ and vice versa: $\vdash_{\textup{EK}} (\textup{M})$ and $\vdash_{EK} (\textup{C})$, logic $\textup{EMC4}$ coincides with the logic $\textup{EK4}=\textup{E4}+\Box(p\rightarrow q)\rightarrow(\Box p\rightarrow\Box q)$. The proof of this fact can be found in \cite{book:Pacuit}.

	\section{Neighborhood semantics} 
	
	In the proof of the main result, we use neighborhood semantics. The idea underlying neighborhood models is implicit in the seminal work by McKinsey and Tarski. Nevertheless, neighborhood semantics were first formally defined by Dana Scott \cite{art:Scott} and Richard Montague \cite{art:Montague} in 1970 independently and in different ways. 
	
	\begin{defi}
		
		A  neighborhood  \textit{frame} is a pair $( W, \bm{\Box} )$, where  $W\neq\varnothing$ and $\bm{\Box}\colon\rho(W)\to\rho(W)$. 
		
	\end{defi}
	
	Usually neighborhood semantics is defined in terms of a function $N: W\rightarrow\rho(\rho(W))$. Given $\bm{\Box}$, we can define $N$ as follows: $N(w)=\{X\mid w\in\Box X\}$. And vice versa: given $N$ we can define function $\bm{\Box}$: $\bm{\Box}X=\{w\mid X\in N(w)\}$. So $\bm{\Box}$ and $N$ are two sides of the same coin. $\bm{\Box}$ is better when we work with iterative axioms like $(4)$. $\bm{\Box}$ allows us to look at the frame `globally' rather than `locally' in case of $N$. 
	
	\begin{defi}
		
		Let $\textup{F}=( W, \bm{\Box} )$ be a neighborhood frame. A  \textit{model} based on $\textup{F}$ is a pair $( \textup{F}, V)$, where $V\colon Var\to\rho(W)$ is a valuation function.
		
	\end{defi}
	
	\begin{defi}\label{def:truth}
		
		Suppose that $M=( \textup{F}, V)$ is a model based on $\textup{F}=( W, \bm{\Box} )$. The truth set of a formula $\varphi$ is defined recursively and denoted by $|\varphi|_{M}$:
		
		\begin{enumerate}[\normalfont 1.]
			\item $|p|_{M}=V(p),~p\in Var$;
			\item $|\neg\varphi|_{M}=W\backslash|\varphi|_{M}$;
			\item $|\varphi\wedge\psi|_{M}=|\varphi|_{M}\cap|\psi|_{M}$;
			\item $|\Box\varphi|_{M}=\bm{\Box}|\varphi|_{M}$.
		\end{enumerate}
	\end{defi}
	
	We say that the formula is true at a point $w\in W$ iff $w\in|\varphi|_{M}$ (in symbols: $M, w\vDash\varphi$).
	
	The proof of the next lemma is straightforward:
	
	\begin{lemma}
		
		Let $M=( \textup{F}, V)$ be a model on $\textup{F}=( W, \bm{\Box} )$, $w\in W$. Then the following is true:

		\begin{enumerate}[\normalfont 1.]
			\item $M, w\vDash p\Leftrightarrow w\in V(p)~(p\in Var)$;
			
			\item  $M, w\vDash\neg\varphi\Leftrightarrow M, w\nvDash\varphi$;
			
			\item $M, w\vDash\varphi\wedge\psi\Leftrightarrow M, w\vDash\varphi$ and $M, w\vDash\psi$.
			
			\item $M, w\vDash \square\varphi\Leftrightarrow w\in\bm{\Box}|\varphi|_{M}$.
		\end{enumerate}
		
	\end{lemma}
	
	A formula $\varphi$ is  \textit{valid} on a frame $\textup{F}=( W, \bm{\Box})$ $(\vDash_{\textup{F}}\varphi)$ iff for every model $M$ on $\textup{F}$, $|\varphi|_{M}=W$. If $\textbf{F}$ is a class of frames and $\textup{L}$ is a logic, we say that $\textbf{F}$  \textit{determines} $\textup{L}$ if the theorems of $\textup{L}$ are exactly the formulas valid on all frames from the class $\textbf{F}$.

	\begin{defi}
		Suppose that $\Gamma\subseteq Fm$ and $\textbf{F}$ is a class of neighborhood frames. A formula $\varphi$ is a  \textit{logical consequence} of $\Gamma$ with respect to $\textbf{F}$ (in symbols, $\Gamma\vDash_{\textbf{F}}\varphi$) if for each model $M$ based on a frame from $\textbf{F}$, for each $w$ in $W$: $M, w\vDash\Gamma$, implies $M,w\vDash\varphi$ (where $M, w\vDash\Gamma$ means $\forall\psi\in\Gamma: M,w\vDash\psi)$.
	\end{defi}
	
	\begin{defi}
		
		A logic $\textup{L}$ is  \textit{strongly complete} with respect to a class of frames $\textbf{F}$ if for each $\Gamma\subseteq Fm$, $\Gamma\vDash_{\textbf{F}}\varphi$ implies $\Gamma\vdash_{L}\varphi$.
		
	\end{defi}
	
	\begin{defi}
		
		A logic $\textup{L}$ has the  \textit{finite model property} (FMP) if there is a class $\textbf{F}$ of finite frames that determines $\textup{L}$.
		
	\end{defi}

	\section{Canonical model}
	
	In this paper we use two the best known techniques of proof completeness: canonical models and filtration. Canonical models sometimes can determine logics. In this case we say that a logic $\textup{L}$ is \textit{canonical}. Filtration is used for establishing the FMP by extracting finite refutation frames from arbitrary frames of a given logic.
	
	In this section we briefly recall the properties of canonical models and then turn to completeness results in the next section. 
	
	The definitions related to canonical models can be found in \cite{book:Pacuit}. We use the abbreviation MCS for `maximally consistent set' and $\textup{L-MCS}$ for maximal $\textup{L}$-consistent set, where $\textup{L}$ is logic.

	\begin{defi}
		
		A  \textit{canonical frame} for a consistent logic $\textup{L}$ is a frame $\textup{F}_{\textup{L}}=( W_{\textup{L}},  \bm{\Box}_{\textup{L}})$ where:
		
		\begin{enumerate}[\normalfont 1.]
			\item $W_\textup{L}=\{\Gamma\mid \Gamma$ is an $\textup{L -- MCS}\}$;
			\item For each $\Gamma\in W_\textup{L}$ and each formula $\varphi\in \textup{L}$:
			$\bm{\Box}_{\textup{L}}|\varphi|_{\textup{L}}=|\Box\varphi|_{\textup{L}}$. (notation: $|\varphi|_{\textup{L}}=\{\Gamma\in W_{\textup{L}}\mid \varphi\in\Gamma\}$).
		\end{enumerate}
		
		Note that $\bm{\Box}_{\textup{L}}X$ for other $X$ is not fixed, so there are many canonical frames for the same logic. 
		
		$M_{\textup{L}}=( \textup{F}_{\textup{L}}, V_{\textup{L}})$ is called a  \textit{canonical model} for a logic $\textup{L}$ if  $\textup{F}_{\textup{L}}$ is the canonical frame for $\textup{L}$ and $V_{\textup{L}}(p)=|p|_{\textup{L}}$ for any propositional variable $p$. 
	\end{defi}
	
	The following lemma and theorem can be also found in \cite{book:Pacuit}.
	
	\begin{lemma}\label{lemma:canlemma}
		Let $\textup{L}$ be a logic and $\varphi,\psi\in\textup{L}$. Then,
		\begin{align*}
			1.&\;|\varphi \wedge \psi|_{\textup{L}}=|\varphi|_{\textup{L}} \cap|\psi|_{\textup{L}} . \\
			2.&\;|\neg\varphi|_{\textup{L}}=M_{\textup{L}}\backslash|\varphi|_{\textup{L}} . \\
			3.&\;|\varphi \vee \psi|_{\textup{L}}=|\varphi|_{\textup{L}} \cup|\psi|_{\textup{L}} \cdot \\
			4.&\;|\varphi|_{\textup{L}} \subseteq|\psi|_{\textup{L}} \text { iff } \vdash_{\textup{L}} \varphi \rightarrow \psi . \\
			5.&\;|\varphi|_{\textup{L}}=|\psi|_{\textup{L}} \text { iff } \vdash_{\textup{L}} \varphi \leftrightarrow \psi .\\
			6.&\;\forall\Gamma:\text{ if }\Gamma \text{ is }\textup{L -- MCS},\varphi \in \Gamma \text{ and }\varphi \rightarrow \psi \in \Gamma \text {, then }\psi\in\Gamma.\\
			7.&\;\forall\Gamma:\text{ if }\Gamma \text{ is }\textup{L -- MCS }\text{and }\vdash_{\textup{L}} \varphi, \text { then } \varphi \in \Gamma.
		\end{align*}
	\end{lemma}

	\begin{theo}\label{theo:truthlemma}
		For any consistent logic $\textup{L}$ and any formula $\varphi$, if $M_{\textup{L}}$ is canonical for $\textup{L}$, then	$|\varphi|_{M}=|\varphi|_{\textup{L}}$.
	\end{theo}
	
	Theorem \ref{theo:truthlemma} implies the following fact: if $\Gamma\nvdash_{\textup{L}}\varphi$ then there is a model $M$ of $\text{L}$ such that $M,\Gamma'\nvDash\varphi$ and $M,\Gamma'\vDash\Gamma$. But still we do not know how a frame of this model looks like. However, to prove the completeness for some logic $\textup{L}$ we have to show that a canonical frame of $\textup{L}$ is in the intended class of frames $\textbf{F}$. 
	
	\begin{defi}
		A logic is called canonical if it is valid on its canonical frame.
	\end{defi}
	
	\begin{prop}
		Every canonical logic is strongly complete.
	\end{prop}

	\section{Completeness}
	
	In this section we prove that $\textup{E4, EMC4, S04}$ are canonical. 
	
	Let us formulate the properties of frames corresponding to the axioms mentioned in the section 2. 
	
	\begin{defi}\label{def:refl}
		
		A frame $\textup{F}=( W,   \bm{\Box})$ is called \emph{reflexive} if for each $X\subseteq W$, $ \bm{\Box} X\subseteq X$.	
	\end{defi}
	
	\begin{defi}\label{def:reg}
		A frame $\textup{F}=( W,   \bm{\Box})$ is called \emph{regular} if for each $X, Y\subseteq W$, $ \bm{\Box} X\cap \bm{\Box} Y\subseteq \bm{\Box}(X\cap Y)$.	
	\end{defi}
	
	The proof of the following lemma is straightforward:
	
	\begin{lemma}
		
		Let $\textup{F}=( W,   \bm{\Box})$  be a regular frame. Then for every $n>0$ and every $X_{1},\ldots, X_{n}\colon \bm{\Box} X_{1}\cap\ldots\cap \bm{\Box}X_{n}\subseteq \bm{\Box}(X_{1}\cap\ldots\cap X_{n}).$
		
	\end{lemma}
	
	\begin{defi}\label{def:mon}
		A frame $\textup{F}=( W,   \bm{\Box})$ is called \emph{monotonic} if for each $X, Y\subseteq W$, $X\subseteq Y\Rightarrow  \bm{\Box} X\subseteq \bm{\Box} Y$.	
	\end{defi}
	
	\begin{defi}\label{def:trans}
		A frame $\textup{F}=( W,   \bm{\Box})$ is called \emph{transitive} if for each $X\subseteq W$, $ \bm{\Box} X\subseteq \bm{\Box} \bm{\Box} X$.
	\end{defi}
	
	Using these definitions we can formulate the following lemma:
	
	\begin{lemma}
		
		Let $\textup{F}=( W, \bm{\Box})$  be a frame. Then:
		\begin{enumerate}[\normalfont (1)]
			\item 	 $\textup{F}\vDash \Box p\rightarrow p$ iff $\textup{F}$ is reflexive;
			
			\item  $\textup{F}\vDash(\Box p\wedge\Box q)\rightarrow\Box(p\wedge q)$ iff $\textup{F}$ is regular;
			
			\item  $\textup{F}\vDash\Box(p\wedge q)\rightarrow(\Box p\wedge\Box q)$ iff $\textup{F}$ is monotonic.
			
			\item $\textup{F}\vDash\Box p\rightarrow\Box\Box p$ iff $\textup{F}$ is transitive.
			
		\end{enumerate}

	\end{lemma}
	
	Let us define the minimal canonical model. It will help us to prove completeness of all three logics from section 2.
	
	\begin{defi}
		The \emph{minimal canonical frame} for $\textup{L}$ is the frame $\textup{F}^{-}_{\textup{L}}=( W_{\textup{L}},  \bm{\Box}_{\textup{L}}^{-})$ such that for every $X\subseteq W$:
		
		$$ \bm{\Box}^{-}_{\textup{L}} X=
		\begin{cases}
			|\Box\varphi|_{\textup{L}}\;, & \text{if}~X=|\varphi|_{\textup{L}}~\text{for some}~\varphi\in Fm;\\
			\varnothing\;, & \text{otherwise}.
		\end{cases} 
		$$
	\end{defi}
	
	Clearly, the definition of $ \bm{\Box}_{\textup{L}}^{-}$ satisfies the conditions of the canonical frame.
	
	\begin{prop}
		Let $\textup{L}$ be a logic. Then $F_{\textup{L}}^{-}$ is canonical.
	\end{prop}
	
	\begin{lemma}\label{can:mon}
		Let $\textup{E4}\subseteq\textup{L}$ and $\textup{F}_{\textup{L}}^{-}$ be the minimal canonical frame for $\textup{L}$. Then $\textup{F}_{\textup{L}}^{-}$ is transitive.
		\begin{proof}
			The proof is straightforward. Let $\Gamma\in \bm{\Box}^{-}_{\textup{L}} X$ for arbitrary $X\subseteq W_{\textup{L}}$. Then there exists $\varphi$ such that $X=|\varphi|_{\textup{L}}$ and $ \bm{\Box}^{-}_{\textup{L}} X=|\Box\varphi|_{\textup{L}}$. Therefore $\Gamma\in |\Box\varphi|_{\textup{L}}$. Since $\Box\varphi\rightarrow\Box\Box\varphi\in\Gamma$, by lemma \ref{lemma:canlemma} we have that $\Gamma\in|\Box\Box\varphi|_{\textup{L}}$. By the definition of the minimal canonical model $\bm{\Box}^{-}_{\textup{L}} \bm{\Box}^{-}_{\textup{L}}X=|\Box\Box\varphi|_{\textup{L}}$. Hence, $\Gamma\in  \bm{\Box}^{-}_{\textup{L}} \bm{\Box}^{-}_{\textup{L}}X$.
		\end{proof}
	\end{lemma} 
	
	To show that $\textup{EMC4}$ and $\textup{S04}$ are complete, we need the following definition:
	
	\begin{defi}\label{def:sup}
		Suppose that $\textup{F}=( W,  \bm{\Box} )$ is a neighborhood frame. The \emph{supplementation} of $\textup{F}$ is $\textup{F}^{\#}=( W,  \bm{\Box}^{\#})$ such that for every $X\subseteq W$:
		$$ \bm{\Box}^{\#}X=\bigcup\{ \bm{\Box} Y\mid Y\subseteq X\}.$$
		
		We say that the model $M^{\#}=( \textup{F}^{\#}, V)$ based on $\textup{F}^{\#}=( W,  \bm{\Box}^{\#})$ is the supplementation of model $M=( \textup{F}, V)$.
	\end{defi}
	
	The proof of the following lemma is straightforward:
	
	\begin{lemma}\label{lemma:mon}
		Let $\textup{F}=( W,  \bm{\Box} )$ be a neigborhood frame. Then its supplementation $\textup{F}^{\#}=( W,  \bm{\Box}^{\#})$ is monotonic.
	\end{lemma}
	
	Therefore we have to prove that the supplementation of the minimal canonical frames for $\textup{S04}$ and $\textup{EMC4}$ are respectively:
	\begin{enumerate}[\normalfont (1)]
		
		\item transitive and reflexive;
		
		\item transitive and regular.
		
	\end{enumerate} 
	
	First let us recall that the minimal canonical models for $\textup{EC}$ and $\textup{ET}$  are regular and reflexive respectively. The proofs of the following lemmas are similar to the proof of Lemma \ref{can:mon}, so we skip them. 
	
	\begin{lemma}\label{lemma:reg}
		Let $\textup{EC}\subseteq\textup{L}$ and $\textup{F}_{\textup{L}}^{-}$ be the minimal canonical frame for $\textup{L}$. Then $\textup{F}_{\textup{L}}^{-}$ is regular. $(\text{Where }\textup{EC=E+C})$
	\end{lemma}
	
	\begin{lemma}\label{Lemma:refl}
		Suppose that $\textup{ET}\subseteq\textup{L}$ and $\textup{F}_{\textup{L}}^{-}$ is the minimal canonical frame for $\textup{L}$. Then $\textup{F}_{\textup{L}}^{-}$ is reflexive. $(\text{Where }\textup{ET=E+T})$
	\end{lemma}
	
	
	Let us prove a general lemma, showing that supplementation preserves reflexivity, transitivity and regularity.
	
	\begin{lemma}\label{lemma:pres}
		Let $\textup{F}=( W, \bm{\Box})$ be a reflexive/transitive/regular frame, then its supplementation $\textup{F}^{\#}=( W,  \bm{\Box}^{\#} )$ is reflexive/transitive/regular.
		\begin{proof}
			First, observe that by lemma \ref{lemma:mon} we have that $Z_{1}\subseteq Z_{2}\Rightarrow \bm{\Box}^{\#}Z_{1}\subseteq \bm{\Box}^{\#}Z_{2}$.
			
			Suppose that $X\subseteq W$. 
			
			\textbf{Case 1: }$\textup{F}$ is reflexive. Suppose that $w\in \bm{\Box}^{\#}X$. Then there exists $Y\subseteq X$ such that $w\in \bm{\Box} Y$. Since $\textup{F}$ is reflexive, $\bm{\Box} Y\subseteq Y$. Therefore, $w\in X$.
			
			\textbf{Case 2: }$\textup{F}$ is transitive. Suppose that $w\in \bm{\Box}^{\#}X$. Then there exists $Y\subseteq X$ such that $w\in \bm{\Box} Y$. By transitivity of $\textup{F}$ we have that $w\in \bm{\Box} \bm{\Box} Y$. Observe that $w\in \bm{\Box}^{\#} \bm{\Box} Y$, since $ \bm{\Box} \bm{\Box} Y\subseteq \bm{\Box}^{\#} \bm{\Box} Y$. By Definition \ref{def:sup} we have that $\bm{\Box} Y\subseteq  \bm{\Box}^{\#}Y$, thus $ \bm{\Box}^{\#} \bm{\Box} Y\subseteq  \bm{\Box}^{\#} \bm{\Box}^{\#}Y$. Therefore, $w\in \bm{\Box}^{\#} \bm{\Box}^{\#}Y$. Since $Y\subseteq X$, we have that $ \bm{\Box}^{\#}Y\subseteq \bm{\Box}^{\#}X$, and also $ \bm{\Box}^{\#} \bm{\Box}^{\#}Y\subseteq \bm{\Box}^{\#} \bm{\Box}^{\#}X$. Hence, $w\in  \bm{\Box}^{\#} \bm{\Box}^{\#}X$. 
			
			\textbf{Case 3: }$\textup{F}$ is regular. Suppose $w\in \bm{\Box}^{\#}X\cap \bm{\Box}^{\#}Y$. Then there are $X'$, $Y'$ such that $w\in\bm{\Box} X', w\in\bm{\Box} Y'\text{ and } X'\subseteq X, Y'\subseteq Y$. By regularity of $\textup{F}$ we obtain that $w\in\bm{\Box}(X'\cap Y')$. Since $X'\cap Y'\subseteq X\cap Y$, we have that $w\in\bm{\Box}^{\#}(X\cap Y)$.
		\end{proof}
	\end{lemma}

	Finally, to prove completeness for $\textup{S04}$ and $\textup{EMC4}$, we need to show that the supplementations of their canonical frames are also canonical. 
	
	\begin{lemma}\label{can:EMC4andS04}
		For $\textup{L}\supseteq\textup{EM}$ the supplementation of its minimal canonical frame $\textup{F}^{-}_{L}$ is canonical.
		\begin{proof}
			Cf.\cite{book:Pacuit}, p.66.
		\end{proof}
	\end{lemma}

	Therefore we obtain a completeness theorem for $\textup{E4}, \textup{EMC4}, \textup{S04}$. Since we proved that these logics are determined by their canonical frames, it follows that they are strongly complete.  
	
	\begin{theo}
		Logics $\textup{E4}, \textup{EMC4}, \textup{S04}$ are strongly complete.
		\begin{proof}
			By lemmas of this section and Theorem \ref{theo:truthlemma}. For details of this proof cf.\cite{book:Pacuit} or \cite{book:BBB}.
		\end{proof}
	\end{theo}
	
	\section{Filtrations}
	
	In this section we will enhance the completeness result and show that every logic from the list above has the FMP. 
	
	First we give the definition of filtrations.Fix a model $M=( F, V)$ on a frame $\textup{F}=( W, \bm{\Box})$ and $\Sigma\subseteq Fm$, which is closed under subformulas. Then define an equivalence relation $\sim_{\Sigma}$ on $W$: $$w\sim_{\Sigma}v:=\forall\varphi\in\Sigma(M, w\vDash\varphi\Leftrightarrow M, v\vDash\varphi).$$
	
	For each $w\in W$ let $\widetilde{w}_{\Sigma}=\{v\mid v\sim_{\Sigma}w\}$ be the equivalence class of $w$. Then $W\slash_{\sim_{\Sigma}}=\{\widetilde{w}_{\Sigma}\mid w\in W\}$, and if $X\subseteq W$, then $\widetilde{X}_{\Sigma}=\{\widetilde{w}_{\Sigma}\mid w\in X\}$. If $\Sigma$ is clear from the context, we may leave out the subscripts. 
	
	\begin{defi}\label{def:fil}
		
		Suppose that $M_{f}=( \textup{F}_{f}, V_{f})$ is a model on a neighborhood frame $\textup{F}_{f}=( W_{f},   \bm{\Box}_{f})$ such that:
		
		\begin{enumerate}[\normalfont 1.]
			
			\item $W_{f}=\widetilde{W}$;
			\item  For every $\varphi$: if $\Box\varphi\in\Sigma$, then $\bm{\Box}_{f}\widetilde{|\varphi|}_{M}=\widetilde{ \bm{\Box}|\varphi|}_{M}$;	
			\item	For every $p$: if $p\in Var,$ then $V_{f}(p)=\widetilde{|p|}_{M}$.
			
		\end{enumerate}
		
		$M_{f}$ is called a \textit{filtration} of the model $M$ through $\Sigma$.
		
	\end{defi}
	
	As in the case of canonical models, $\bm{\Box}_{f}$ can vary on sets $X\subseteq W_{f}$ that are not of the form $\widetilde{|\varphi|}_{M}$. Thus, there are many different filtrations. We shall see the examples below.
	
	We have the following important theorem for filtrations:
	
	\begin{theo}
		Let $M=( \textup{F}, V)$ be a model and $M_{f}=( \textup{F}_{f}, V_{f} )$ be a filtration of $M$ through $\Sigma$. Then, for each $\varphi\in\Sigma$,
		$$|\varphi|_{M_{f}}=\widetilde{|\varphi|}_{M}.$$ 
		\begin{proof}
			Cf.\cite{book:Pacuit}, p.70.
		\end{proof}
	\end{theo}
	
	The following technical lemma follows:
	
	\begin{lemma}\label{lemma:tech1}
		Let $M_{f}=( \textup{F}_{f}, V_{f} )$ be a filtration of model $M=( \textup{F}, V )$ through $\Sigma$. For every $\varphi, \psi\in\Sigma$, $|\varphi|_{M}=|\psi|_{M}\Leftrightarrow\widetilde{|\varphi|}_{M}=\widetilde{|\psi|}_{M}.$
		\begin{proof}
			Follows from the definition \ref{def:fil}.
		\end{proof}
	\end{lemma}
	
	The definitions \ref{def:refl}, \ref{def:reg}, \ref{def:mon}, \ref{def:trans}, \ref{def:sup} are naturally extrapolated on models.
	
	Let us show that $\textup{E4}$ has the FMP. To prove this we introduce the concept of minimal filtration. This filtration is the smallest filtration in the lattice of filtrations for every model. However, we are not going to prove this fact here. 
	
	\begin{defi}\label{def:min.fil}
		
		Suppose that $M=( \textup{F}, V)$ is a model on a neighborhood frame $\textup{F}=( W,   \bm{\Box})$. The \textit{minimal filtration} of $M$ through set $\Sigma$ is model $M^{-}_{f}=( \textup{F}^{-}_{f}, V_{f} )$ such that $\textup{F}^{-}_{f}=( W_{f},  \bm{\Box}_{f}^{-} )$ and for each $X\subseteq W_{f}$:
		
		$$ \bm{\Box}^{-}_{f}X=
		\begin{cases}
			\widetilde{|\Box\varphi|}_{M},& \text{if}~X=\widetilde{|\varphi|}_{M}\text{ for some formula }\Box\varphi\in\Sigma;\\
			\varnothing,& \text{otherwise}.
		\end{cases}
		$$
	\end{defi}
	
	We have to check that minimal filtration is well-defined. But first note that the following technical proposition holds:
	
	\begin{prop}
		
		Let $M$ be a model and $M^{-}_{f}$ be minimal filtration of $M$ trough $\Sigma$. Then for every formulas $\Box\varphi, \Box\psi\in\Sigma:$ if $\widetilde{|\psi|}_{M}=\widetilde{|\varphi|}_{M}$, then $\bm{\Box}_{f}^{-}\widetilde{|\psi|}_{M} =\bm{\Box}_{f}^{-}\widetilde{|\varphi|}_{M}$.
		
		\begin{proof}	
			It is sufficient to show that $\widetilde{|\Box\varphi|}_{M}=\widetilde{|\Box\psi|}_{M}$. Let us prove only one inclusion, the other one is completely analogous.
			
			Assume that $\widetilde{w}\in \widetilde{|\Box\varphi|}_{M}$. Then there is $v\sim w$, s.t. $M, v\vDash\Box\varphi$. Then $v\in \bm{\Box}|\varphi|_{M}$. Since $\widetilde{|\varphi|}_{M}=\widetilde{|\psi|}_{M}$, then we also have that $|\varphi|_{M}=|\psi|_{M}$. Hence, $\bm{\Box}|\varphi|_{M}=\bm{\Box}|\psi|_{M}$. So we have that $v\in  \bm{\Box}|\psi|_{M}$. Thus, $v\in |\Box\psi|_{M}$ and therefore $\widetilde{v}\in\widetilde{|\Box\psi|}_{M}$. Thus, $\widetilde{w}\in\widetilde{|\Box\psi|}_{M}$.
		\end{proof}
		
	\end{prop}
	
	Observe, that the minimal filtration is indeed a filtration according to Definition \ref{def:fil}.
	
	Now we introduce the key definition for construction of filtration for transitive models:
	
	\begin{defi}\label{def:closure}
		Let $\textup{F}=( W,   \bm{\Box})$ be a neighborhood frame. The \textit{closure of the function} $ \bm{\Box}$ is $\widehat{ \bm{\Box}}\colon\rho(W)\to\rho(W)$ such that for every $X\subseteq W$:
		
		$$\widehat{ \bm{\Box}}X=
		\begin{cases}
			X,& \text{if}~X= \bm{\Box} Y~\text{for some}~Y\subseteq W;\\
			\varnothing,& \text{otherwise}.
		\end{cases}
		$$
		
	\end{defi}
	
	Using Definition \ref{def:closure} we can now define filtrations for transitive models.
	
	\begin{defi}\label{def:trans.fil}
		
		Suppose that $M=( \textup{F}, V)$ is a model and $M^{-}_{f}=( \textup{F}^{-}_{f}, V_{f} )$  is the minimal filtration of $M$ through $\Sigma$. \textit{The transitive filtration} of $M$ is the model $M^{T}_{f}=( \textup{F}^{T}_{f}, V_{f})$ based on the frame $\textup{F}^{T}_{f}=( W_{f},  \bm{\Box}_{f}^{T} )$ such that for every $X\subseteq W_{f}$:
		
		$$ \bm{\Box}_{f}^{T}X= \bm{\Box}_{f}^{-}X\cup\widehat{ \bm{\Box}^{-}_{f}}X.$$
		
	\end{defi}
	
	First let us prove that $M^{T}_{f}$ is a filtration according to Definition \ref{def:fil}.
	
	For this purpose we need the following technical lemma:
	
	\begin{lemma}\label{lem:tech}
		Let $M_{f}=( \textup{F}_{f}, V_{f} )$ be a filtration of model $M=( \textup{F}, V )$ through $\Sigma$. If $\widetilde{|\varphi|}_{M}= \bm{\Box}_{f}\widetilde{|\psi|}_{M}$ and $\varphi, \Box\psi\in\Sigma$, then $\widetilde{|\Box\varphi|}_{M}=\widetilde{|\Box\Box\psi|}_{M}$.
		\begin{proof}
			By the Lemma \ref{lemma:tech1} we have that $|\varphi|_{M}=|\Box\psi|_{M}$, so $\bm{\Box}|\varphi|_{M}=\bm{\Box}|\Box\psi|_{M}$. So, by definition of truth in a model: $|\Box\varphi|_{M}=|\Box\Box\psi|_{M}$. Therefore, $\widetilde{|\Box\varphi|}_{M}=\widetilde{|\Box\Box\psi|}_{M}$.
			
		\end{proof}
	\end{lemma}
	
	\begin{lemma}
		Let $M_{f}^{-}=( \textup{F}_{f}^{-}, V_{f} )$ be the minimal filtration of a transitive model $M=( \textup{F}, V )$ through $\Sigma$. Then $M_{f}^{T}=( \textup{F}_{f}^{T}, V_{f} )$ is a filtration of the model $M$ through $\Sigma$.
		\begin{proof}
			
			We must show that $\bm{\Box}_{f}^{T}\widetilde{|\varphi|}_{M}=\widetilde{|\Box\varphi|}_{M}$ for every formula $\Box\varphi\in\Sigma$.
			
			Let $\Box\varphi\in\Sigma$.\\	
			$(\subseteq)$ Suppose that $\widetilde{w}\in \bm{\Box}^{T}_{f}\widetilde{|\varphi|}_{M}$. Then $\widetilde{w}\in \bm{\Box}^{-}_{f}\widetilde{|\varphi|}_{M}$ or $\widetilde{w}\in\widehat{ \bm{\Box}^{-}_{f}}\widetilde{|\varphi|}_{M}$. In the first case $\widetilde{w}\in\widetilde{|\Box\varphi|}_{M}$, since $M^{-}_{f}$ is a filtration. In the second case we have that for some $Y$, $\widetilde{|\varphi|}_{M}= \bm{\Box}^{-}_{f}Y$ and $\widetilde{w}\in \widetilde{|\varphi|}_{M}$. Therefore, there is a formula $\Box\psi\in\Sigma$ such that $Y=\widetilde{|\psi|}_{M}$ and $ \bm{\Box}_{f}^{-}Y=\widetilde{|\Box\psi|}_{M}$. Hence, $\widetilde{w}\in\widetilde{|\Box\psi|}_{M}$. Using transitivity of $M$, we obtain that $\widetilde{w}\in \widetilde{|\Box\Box\psi|}_{M}$. Since  $\widetilde{|\varphi|}_{M}= \widetilde{|\Box\psi|}_{M}$ and $\varphi, \Box\psi\in\Sigma$, using Lemma \ref{lem:tech} we obtain that $\widetilde{|\Box\varphi|}_{M}=\widetilde{|\Box\Box\psi|}_{M}$.  Hence, $\widetilde{w}\in\widetilde{|\Box\varphi|}_{M}$.\\
			$(\supseteq)$ If $\widetilde{w}\in\widetilde{|\Box\varphi|}_{M}$, then $\widetilde{w}\in \bm{\Box}_{f}^{-}\widetilde{|\varphi|}_{M}$. Therefore $\widetilde{w}\in \bm{\Box}_{f}^{T}\widetilde{|\varphi|}_{M}$.
		\end{proof}
	\end{lemma}
	
	Now we will show that applying the transitive filtration to a transitive model preserves transitivity.
	
	\begin{lemma}
		Let $M=( \textup{F}, V )$ be a transitive model. Then $M_{f}^{T}$ is transitive.
		
		\begin{proof}
			
			We must show that for every $X$, $ \bm{\Box}^{T}_{f}X\subseteq \bm{\Box}^{T}_{f} \bm{\Box}^{T}_{f}X$.
			
			By Definition \ref{def:trans.fil}: $ \bm{\Box}^{T}_{f} \bm{\Box}^{T}_{f}X= \bm{\Box}^{T}( \bm{\Box}^{-}_{f} X\cup\widehat{ \bm{\Box}^{-}_{f}}X )= \bm{\Box}^{-}_{f}( \bm{\Box}^{-}_{f} X\cup\widehat{ \bm{\Box}^{-}_{f}}X )\cup\widehat{ \bm{\Box}_{f}^{-}}( \bm{\Box}^{-}_{f} X\cup\widehat{ \bm{\Box}^{-}_{f}}X )$.
			
			There are two cases: $\bm{\Box}^{-}_{f}X=\varnothing$ or $\bm{\Box}^{-}_{f}X\neq\varnothing$.
			
			\textbf{Case 1:} If $ \bm{\Box}^{-}_{f}X=\varnothing$ and $\widetilde{w}\in  \bm{\Box}^{T}_{f}X$, then $\widetilde{w}\in\widehat{ \bm{\Box}^{-}_{f}}X$. Therefore, $\widehat{\bm{\Box}_{f}^{-}}X=X$. Then $ \bm{\Box}^{T}_{f} \bm{\Box}^{T}_{f}X= X$. Hence, $\widetilde{w}\in  \bm{\Box}^{T}_{f} \bm{\Box}^{T}_{f}X$.
			
			\textbf{Case 2:} Suppose that $\bm{\Box}^{-}_{f}X\neq\varnothing$. We will show that in this case $ \bm{\Box}^{-}_{f}X=\bm{\Box}^{T}_{f}X$. 
			
			$(\subseteq)$ This inclusion follows from the definition of $ \bm{\Box}^{T}_{f}$.
			
			$(\supseteq)$ Suppose that $\widetilde{w}\in \bm{\Box}^{T}_{f}X$. Then $ \widetilde{w}\in \bm{\Box}^{-}_{f}X$ or $\widetilde{w}\in\widehat{ \bm{\Box}^{-}_{f}}X$. If $\widetilde{w}\in \bm{\Box}^{-}_{f}X$, then  the claim is proved. If $\widetilde{w}\in\widehat{ \bm{\Box}^{-}_{f}}X$, then there is $Y$ such that $X= \bm{\Box}^{-}_{f}Y$. Therefore, there exists $\Box\psi\in\Sigma$ such that $Y=\widetilde{|\psi|}_{M}$ and $\bm{\Box}^{-}_{f}Y=\widetilde{|\Box\psi|}_{M}$. Then, $ \bm{\Box}^{-}_{f}X= \bm{\Box}^{-}_{f} \bm{\Box}^{-}_{f}\widetilde{|\psi|}_{M}$. Since $\widetilde{w}\in \bm{\Box}^{-}_{f}\widetilde{|\psi|}_{M}$, then we have that $\widetilde{w}\in\widetilde{|\Box\psi|}_{M}$. By transitivity of $M$, we have that $\widetilde{w}\in\widetilde{|\Box\Box\psi|}_{M}$. 
			
			Since $ \bm{\Box}^{-}_{f}X\neq\varnothing$, we have that there exists $\Box\varphi\in\Sigma$ such that $X=\widetilde{|\varphi|}_{M}$ and $\bm{\Box}^{-}_{f}X=\widetilde{|\Box\varphi|}_{M}$. Using the fact that $\widetilde{|\Box\psi|}_{M}=\widetilde{|\varphi|}_{M}$ and $\Box\varphi, \Box\psi\in\Sigma$, we obtain by Lemma \ref{lem:tech} that $\widetilde{|\Box\varphi|}_{M}=\widetilde{|\Box\Box\psi|}_M$. Therefore, $\widetilde{w}\in \widetilde{|\Box\varphi|}_{M}$. Hence, $\widetilde{w}\in \bm{\Box}^{-}_{f}\widetilde{|\varphi|}_{M}$, and so,  $\widetilde{w}\in \bm{\Box}^{-}_{f}X$.
			
			Using the proved equality $ \bm{\Box}^{-}_{f}X= \bm{\Box}^{T}_{f}X$ the proof of the inclusion $\bm{\Box}^{T}_{f}X\subseteq \bm{\Box}^{T}_{f} \bm{\Box}^{T}_{f}X$  is straightforward: since $ \bm{\Box}^{-}_{f}X\subseteq \bm{\Box}^{-}_{f} \bm{\Box}^{-}_{f}X\cup \bm{\Box}^{-}_{f}X= \bm{\Box}^{-}_{f} \bm{\Box}^{-}_{f}X\cup\widehat{ \bm{\Box}^{-}_{f}} \bm{\Box}_{f}^{-}X= \bm{\Box}^{T}_{f} \bm{\Box}^{-}_{f}X$, we have $ \bm{\Box}^{T}_{f}X\subseteq \bm{\Box}^{T}_{f} \bm{\Box}^{T}_{f}X$. Note that here we use that $\bm{\Box}^{T}_{f}$ is a function.

			Therefore in both cases $ \bm{\Box}^{T}_{f}X\subseteq \bm{\Box}^{T}_{f} \bm{\Box}^{T}_{f}X$, so $M^{T}_{f}$ is transitive.
		\end{proof}
	\end{lemma}
	
	We can now move on to the logic $\textup{S04}$, the logic of closure operators as stated in the introduction. We will show that it also has the FMP.
	
	\begin{lemma}
		Let $M=( \textup{F}, V )$ be a monotonic and transitive model and $M_{f}^{T}$ be the transitive filtration of $M$ through $\Sigma$. Then its supplementation  $M_{f}^{T\#}=( \textup{F}_{f}^{T\#}, V_{f} )$ is a filtration of $M$ through $\Sigma$.
		\begin{proof}
			
			We must show that $ \bm{\Box}_{f}^{T\#}\widetilde{|\varphi|}_{M}=\widetilde{|\Box\varphi|}_{M}$ for every formula  $\Box\varphi\in\Sigma$.
			
			Let $\Box\varphi\in\Sigma$.\\
			$(\subseteq)$ Suppose that $\widetilde{w}\in  \bm{\Box}_{f}^{T\#}\widetilde{|\varphi|}_{M}$. Then there exists $Y\subseteq \widetilde{|\varphi|}_{M}$ such that $\widetilde{w}\in \bm{\Box}_{f}^{T}Y$. Then $\widetilde{w}\in \bm{\Box}^{-}_{f}Y$ or $\widetilde{w}\in\widehat{ \bm{\Box}^{-}_{f}}Y$.
			
			If $\widetilde{w}\in \bm{\Box}^{-}_{f}Y$, then there is $\Box\psi\in\Sigma$ such that $Y=\widetilde{|\psi|}_{M}$ and $ \bm{\Box}^{-}_{f}Y=\widetilde{|\Box\psi|}_{M}$. Since $\widetilde{|\psi|}_{M}\subseteq\widetilde{|\varphi|}_{M}$ and $\Box\psi, \Box\varphi\in\Sigma$, then by definition of our equivalence relation $\sim_\{\Sigma\}$ we have that $|\psi|_{M}\subseteq |\varphi|_{M}$. Monotonicity of $\bm{\Box}$ in model $M$ gives us $\bm{\Box}|\psi|_{M}\subseteq \bm{\Box}|\varphi|_{M}$ Hence, $|\Box\psi|_{M}\subseteq|\Box\varphi|_{M}$. Thus, $\widetilde{|\Box\psi|}_{M}\subseteq\widetilde{|\Box\varphi|}_{M}$ and  $\widetilde{w}\in \widetilde{|\Box\varphi|}_{M}$. 
			
			If $\widetilde{w}\in\widehat{ \bm{\Box}^{-}_{f}}Y$, then there exists $U$ such that $Y= \bm{\Box}^{-}_{f}U$. Thus, there is $\Box\gamma\in\Sigma$ such that $U=\widetilde{|\gamma|}_{M}$ and $ \bm{\Box}^{-}_{f}U=\widetilde{|\Box\gamma|}_{M}$. Since $M$ is transitive, we have that $\widetilde{|\Box\gamma|}_{M}\subseteq\widetilde{|\Box\Box\gamma|}_{M}$, and so, $\widetilde{w}\in\widetilde{|\Box\Box\gamma|}_{M}$. Therefore, using the fact that $\widetilde{|\Box\gamma|}_{M}=Y\subseteq\widetilde{|\varphi|}_{M}$, $\Box\gamma,\Box\varphi\in\Sigma$ and by monotonicity of $M$ we obtain that $\widetilde{|\Box\Box\gamma|}_{M}\subseteq\widetilde{|\Box\varphi|}_{M}$. 
			Therefore, $\widetilde{w}\in\widetilde{|\Box\varphi|}_{M}$.\\
			$(\supseteq)$ If $\widetilde{w}\in\widetilde{|\Box\varphi|}_{M}$. Then $\widetilde{w}\in \bm{\Box}^{T}_{f}\widetilde{|\varphi|}_{M}$. Consequently, $\widetilde{w}\in \bm{\Box}^{T\#}_{f}\widetilde{|\varphi|}_{M}$. 
		\end{proof}
	\end{lemma} 
	
	Now we shall prove that the transitive filtration of a reflexive frame is reflexive. 
	
	\begin{lemma}
		Let $M=( \textup{F}, V )$ be a reflexive model. Then its transitive filtration $M^{T}_{f}$ through $\Sigma$ is reflexive.
		\begin{proof}
			
			We must show that $\bm{\Box}^{T}_{f}X\subseteq X$. Suppose that	$\widetilde{w}\in \bm{\Box}^{T}_{f}X$ for arbitrary $X\subseteq W_{f}$. Then $\widetilde{w}\in \bm{\Box}^{-}_{f}X$ or $\widetilde{w}\in\widehat{ \bm{\Box}^{-}_{f}}X$. Let us consider the first case when $\widetilde{w}\in \bm{\Box}^{-}_{f}X$. If so, then there exists $\Box\varphi\in\Sigma$ such that $ \bm{\Box}^{-}_{f}X=\widetilde{|\Box\varphi|}_{M}$ and $X=\widetilde{|\varphi|}_{M}$.  Then using that $\Box\varphi\in\Sigma$ and that $M$ is reflexive, we have $\widetilde{w}\in\widetilde{|\varphi|}_{M}$, since $M$ is reflexive. Therefore, $\widetilde{w}\in X$. 
			
			In the second case when $\widetilde{w}\in\widehat{ \bm{\Box}^{-}_{f}}X$, we have that $\widehat{ \bm{\Box}^{-}_{f}}X=X$, thus $\widetilde{w}\in X$.	
		\end{proof}
	\end{lemma}
	
	To finish the proof of the FMP for $\textup{S04}$ it remains to re that supplementation preserves reflexivity and transitivity.

	To prove the FMP for $\textup{EMC4}$ we will need the definition of the intersection closure. 
	
	\begin{defi}
		Suppose that $\textup{F}=( W,  \bm{\Box} )$ is a neighborhood frame. The \textit{intersection closure} of $\textup{F}$ is the frame $\textup{F}^{*}=( W,  \bm{\Box}^{*} )$ such that for every $X\subseteq W$:
		$$ \bm{\Box}^{*}X=\bigcup\{ \bm{\Box} X_{1}\cap\ldots\cap \bm{\Box} X_{n}\mid X=X_{1}\cap\ldots\cap X_{n}\}.$$
	\end{defi}
	
	
	The proof of the next theorem is straightforward.
	
	\begin{theo}
		
		Let $\textup{F}=( W,  \bm{\Box} )$ be a neighborhood frame. Then its intersection closure $\textup{F}^{*}=( W,  \bm{\Box}^{*} )$ is regular.
		
	\end{theo}
	
	To show that the regular and monotonic models admit filtration we will use the construction presented in Brian Chellas' book \cite{book:Chellas}. 
	
	\begin{lemma}
		
		Suppose that $\textup{F}=( W,  \bm{\Box} )$ is a neighborhood frame. Then $\textup{F}^{*\#}=\textup{F}^{\#*}$.
		\begin{proof}
			We must show that for arbitrary $X\subseteq W$, $ \bm{\Box}^{*\#}X= \bm{\Box}^{\#*}X$.
			
			$(\subseteq)$ Suppose that $w\in \bm{\Box}^{*\#}X$. It is sufficient to show that $\exists U_{1},\ldots, U_{n}: U_{1}\cap\ldots\cap U_{n}=X\text{ and }w\in \bm{\Box}^{\#}U_{1}\cap\ldots\cap \bm{\Box}^{\#}U_{n}$. Since $w\in \bm{\Box}^{*\#}X$ there exists $Y\subseteq X$ such that $w\in \bm{\Box}^{*}Y$. Therefore there are $Y_{1},\ldots, Y_{n}$ such that $Y=Y_{1}\cap\ldots\cap Y_{n}$ and $w\in \bm{\Box} Y_{1}\cap\ldots\cap \bm{\Box} Y_{n}$. Note that $X=Y\cup(X\setminus Y)=(Y_{1}\cap\ldots\cap Y_{n})\cup(X\setminus Y)=(Y_{1}\cup(X\setminus Y))\cap\ldots\cap(Y_{n}\cup(X\setminus Y))$. Since $Y_{i}\subseteq (Y_{i}\cup(X\setminus Y))$, we have that $ \bm{\Box} Y_{i}\subseteq \bm{\Box}^{\#}(Y_{i}\cup(X\setminus Y))$. Then we can take $U_{i}=Y_{i}\cup (X\setminus Y)$. Thus, $w\in  \bm{\Box}^{\#*}X$.
			
			$(\supseteq)$ If $w\in  \bm{\Box}^{\#*}X$, then there are $X_{1},\ldots, X_{n}$ such that  $X=X_{1}\cap\ldots\cap X_{n}$ and $w\in \bm{\Box}^{\#} X_{1}\cap\ldots\cap \bm{\Box}^{\#} X_{m}$. Therefore for every $X_{i}$ there exists $Y_{i}$ such that $Y_{i}\subseteq X_{i}$ and $w\in \bm{\Box} Y_{i}$. Therefore $w\in \bm{\Box}^{*}Y$ for $Y=Y_{1}\cap\ldots\cap Y_{m}$. Since $Y\subseteq X$ it follows that $w\in  \bm{\Box}^{*\#}X$.
		\end{proof} 
	\end{lemma}
	
	Now we can define the composition of these operations in the following way:
	
	\begin{defi}
		Suppose that $\textup{F}=( W,  \bm{\Box} )$ is a neighborhood frame. The \textit{rm-closure} of $\textup{F}$ is the frame $\textup{F}^{\bullet}=( W,  \bm{\Box}^{\bullet})$ such that for every $X\subseteq W$:  $$ \bm{\Box}^{\bullet}X= \bm{\Box}^{*\#}X= \bm{\Box}^{\#*}X.$$
	\end{defi}
	
	Now we will show that the rm-closure preserves transitivity.
	
	\begin{lemma}
		Let $\textup{F}=( W,  \bm{\Box} )$ be a monotonic and transitive neighborhood frame. Then its rm-closure $\textup{F}^{\bullet}$ is transitive. 
		\begin{proof}
			We have to prove that for an arbitrary $X\subseteq W$, $ \bm{\Box}^{\bullet}X\subseteq \bm{\Box}^{\bullet} \bm{\Box}^{\bullet}X$. 
			Suppose that $w\in  \bm{\Box}^{\bullet}X$. Then there exists $Y\subseteq X$ such that $w\in \bm{\Box}^{*} Y$. Therefore there are $Y_{1},\ldots, Y_{n}$ such that  $Y=Y_{1}\cap\ldots\cap Y_{n}$ and $w\in \bm{\Box} Y_{1}\cap\ldots\cap \bm{\Box} Y_{n}$. Since $M$ is transitive, we have $w\in \bm{\Box} \bm{\Box} Y_{1}\cap\ldots\cap \bm{\Box} \bm{\Box} Y_{n}$. Therefore $w\in \bm{\Box}^{*}( \bm{\Box} Y_{1}\cap\ldots\cap\bm{\Box} Y_{n})$, and so, $w\in \bm{\Box}^{\bullet}( \bm{\Box} Y_{1}\cap\ldots \cap\bm{\Box} Y_{n})$. Since we have inclusions $ \bm{\Box} Y_{1}\cap\ldots\cap \bm{\Box} Y_{n}\subseteq \bm{\Box}^{*}Y\subseteq \bm{\Box}^{\bullet}X$ and $ \textup{F}^{\bullet}$ is monotonic, it follows that $ \bm{\Box}^{\bullet}( \bm{\Box} Y_{1}\cap\ldots \bm{\Box} Y_{n})\subseteq \bm{\Box}^{\bullet} \bm{\Box}^{\bullet}X$. Hence, $w\in \bm{\Box}^{\bullet} \bm{\Box}^{\bullet}X$.
		\end{proof}
	\end{lemma}
	
	\begin{lemma}
		
		Let $M=( \textup{F}, V )$ be a transitive, monotonic and regular model and let $M^{T}_{f}$ be its transitive filtration through $\Sigma$. Then $M^{T\bullet}_{f}$ is a filtration of $M$ through $\Sigma$.
		
		\begin{proof}
			
			We must show that $ \bm{\Box}_{f}^{T\bullet}\widetilde{|\varphi|}_{M}=\widetilde{|\Box\varphi|}_{M}$ for every formula $\Box\varphi\in\Sigma$.	
			$(\subseteq)$ Assume that $\widetilde{w}\in \bm{\Box}_{f}^{T\bullet}\widetilde{|\varphi|}_{M}$. Then there exists $Y\subseteq\widetilde{|\varphi|}_{M}$ such that $\widetilde{w}\in \bm{\Box}^{T*}_{f}Y$. Then there are $Y_{1},\ldots, Y_{n}$ such that  $Y=Y_{1}\cap\ldots\cap Y_{n}$ and $\widetilde{w}\in \bm{\Box}^{T}_{f} Y_{1}\cap\ldots\cap \bm{\Box}^{T}_{f} Y_{n}$.Therefore, for every $i\in\{1,\ldots, n\}$, $\widetilde{w}\in \bm{\Box}^{-}_{f} Y_{i}$ or $\widetilde{w}\in\widehat{ \bm{\Box}^{-}_{f}} Y_{i}$. Let us divide the set $\{Y_{1},\ldots, Y_{n}\}$ into two parts. The first part consists of $Y_{i}$ such that $\widetilde{w}\in \bm{\Box}^{-}_{f}Y_{i}$ (we denote these $Y_{i}$ by $V_{j}$, where $j\in\{1,\ldots,l\}$). The second part consists of $Y_{i}$ such that $\widetilde{w}\in\widehat{ \bm{\Box}^{-}_{f}} Y_{i}$ (we denote these $Y_{i}$ by $U_{k}$, where $k\in\{l+1,\ldots,n\}$). For every $V_{j}$ there exists $\Box\psi_{j}\in\Sigma$ such that $V_{j}=\widetilde{|\psi_{j}|}_{M}$ and $ \bm{\Box}^{-}_{f}V_{j}=\widetilde{|\Box\psi_{j}|}_{M}$. On the other hand, for every $U_{k}$ there exists $Z_{k}$ such that $U_{k}= \bm{\Box}_{f}^{-}Z_{k}$. Hence, for every $Z_{k}$, there exists $\Box\gamma_{k}\in\Sigma$ such that $Z_{k}=\widetilde{|\gamma_{k}|}_{M}$ and $ \bm{\Box}_{f}^{-}Z_{k}=\widetilde{|\Box\gamma_{k}|}_{M}$. Since $Y=V_{1}\cap\ldots\cap V_{l}\cap U_{l+1}\cap\ldots\cap U_{n}$, we obtain $Y=\bigcap\limits_{j=1}^{l}\widetilde{|\psi_{j}|}_{M}\cap\bigcap\limits_{k=l+1}^{n}\widetilde{|\Box\gamma_{k}|}_{M}$. Since $Y\subseteq \widetilde{|\varphi|}_{M}$, it is clear that $\bigcap\limits_{j=1}^{l}|\psi_{j}|_{M}\cap\bigcap\limits_{k=l+1}^{n}|\Box\gamma_{k}|_{M}\subseteq|\varphi|_{M}$ (note that $\Box\varphi, \Box\gamma_{k}, \Box\psi_{j}\in\Sigma$ for every $k,j$). By monotonicity of $\textup{F}$, $ \bm{\Box}(\bigcap\limits_{j=1}^{l}|\psi_{j}|_{M}\cap\bigcap\limits_{k=l+1}^{n}|\Box\gamma_{k}|_{M})\subseteq \bm{\Box}|\varphi|_{M}$. 
			
			Now we will show that $w\in  \bm{\Box}(\bigcap\limits_{j=1}^{l}|\psi_{j}|_{M}\cap\bigcap\limits_{k=l+1}^{n}|\Box\gamma_{k}|_{M})$. 
			
			Since  $\widetilde{w}\in \widetilde{|\Box\gamma_{k}|}_{M}$ and $\widetilde{w}\in\widetilde{|\Box\psi_{j}|}_{M}$, we have $w\in|\Box\gamma_{k}|_{M}$ and $w\in|\Box\psi_{j}|_{M}$ for every $k, j$. Since $\textup{F}$ is transitive, we obtain $w\in|\Box\Box\gamma_{k}|_{M}$ for every $k$. Therefore $\widetilde{w}\in\bigcap\limits_{j=1}^{l}|\Box\psi_{j}|_{M}\cap\bigcap\limits_{k=l+1}^{n}|\Box\Box\gamma_{k}|_{M}$ for every $j,k$. By Definition \ref{def:truth}, $w\in\bigcap\limits_{j=1}^{l} \bm{\Box}|\psi_{j}|_{M}\cap\bigcap\limits_{k=l+1}^{n} \bm{\Box}|\Box\gamma_{k}|_{M}$. Hence, by regularity of $\textup{F}$, $w\in  \bm{\Box}(\bigcap\limits_{j=1}^{l}|\psi_{j}|_{M}\cap\bigcap\limits_{k=l+1}^{n}|\Box\gamma_{k}|_{M})$.
			
			Therefore, $w\in \bm{\Box}|\varphi|_{M}=|\Box\varphi|_{M}$, so $\widetilde{w}\in\widetilde{|\Box\varphi|}_{M}$.
			
			$(\supseteq)$ If $\widetilde{w}\in\widetilde{|\Box\varphi|}_{M}$, then $\widetilde{w}\in \bm{\Box}^{T}_{f}\widetilde{|\varphi|}_{M}$ and therefore $\widetilde{w}\in \bm{\Box}^{T\bullet}_{f}\widetilde{|\varphi|}_{M}$.
		\end{proof}
		
	\end{lemma}
	
	Then we have the following theorem:
	
	\begin{theo}\label{theo:main}
		
		The	logics $\textup{E4}, \textup{S04}, \textup{EMC4}$ have the FMP. 
		
	\end{theo}
	
	Using Theorem \ref{theo:main} we can extend the result and establish that every logic obtained by adding a variable free formula $\psi$ to $\textup{E4}, \textup{S04}$ or $\textup{EMC4}$ has the FMP.
	
	\begin{coro}
		
		For every variable free formula $\psi$, the logics $\textup{E4}+\psi$,  $\textup{S04}+\psi$, $\textup{EMC4}+\psi$ admit filtration, and so, have the the FMP.
		
	\end{coro}
	
	\section{Conclusion and Future Work}
	
	We have shown that some non-normal modal logics with transitivity have the FMP. The question about the FMP of logic $\textup{EMC4}$ was posed in the paper \cite{art:Surendonk} by T.J.Surendonk. Another interesting logic with FMP is $\textup{S04}$, since its axioms reflect the behavior of closure operators. However, the main construction for the proof of the FMP of logics with transitivity was developed while working with the logic $\textup{E4}$. Thus, the key concept is the definition of \textit{transitive filtration}. It allows us to build a filtration, which is transitive. Note that applying transitive filtration to some arbitrary model might not give us a transitive model. Thus, the first open problem that we would like to stress is the uniform definition of transitive closure for neighborhood models. One of the possible directions for solving this problem could be to examine the connection between neighborhood models and relational models (cf.\cite{book:Pacuit}). Since the question of transitive closure for the relational models is studied extensively, it might give new insights with the following procedures: 1) transform neighborhood model into relational model, 2) apply a transitive closure and 3) go back to neighborhood model. 
	
	As stated in the introduction, Lewis's result about (weak) completeness and the FMP of the class of non-iterative modal logics with finitely many axioms, contains an error. However, Surendonk proved the strong completeness property for class of non-iterative modal logics via algebraic semantics. Since in the proof of Surendonk there is no analogy to Lewis' error, it might be possible to fix Lewis' proof. Having the completeness and the FMP result for the class of non-iterative modal logics, we can pose another open question. Suppose that $\textup{L}$ is a non-iterative monomodal logic, is $\textup{L} + 4$ complete? Does it have the FMP? What changes if $\textup{L}$ is poly-modal? Even if the answer for the first question is positive, the answer for the second question is not so obvious. As we did not find a filtration for logic \textup{EC4}, it might be worth considering this logic at first.
	
	\appendix
	
	\section{}\label{appendix.A}
	
	We are going to examine the proof of Lewis until the error. It should be mentioned that Lewis' notation is quite outdated. Also, he works with modal operator of arity $n$. However, restricting the proof to the case with $\Box$ does not change the essence of the proof, so we are going to work with it. One can find the same proof with modern notation in the textbook \cite{book:Pacuit}. The error in this case is in the proposition 2.7 on page 74. 
	
	Lewis strategy is to build a model and using it prove completeness and the FMP for some logic $\textup{L}$, axiomatized by finitely many non-iterative formulas. The construction is the following:
	\begin{enumerate}
		\item Take a formula $\varphi$ which we want to refute in our future model;\\
		\item Having a set $\Sigma\subseteq Sub(\varphi)$, define $\varphi$ -- description $D$ as follows: $D=\Sigma\cup \{\neg\psi\mid \psi\in Sub(\varphi)\backslash\Sigma\}$;\\
		\item For each L -- consistent $\varphi$ -- descriprion $D$ choose an L -- MCS $\Gamma_{D}$ such that $D\subseteq \Gamma_{D}$;
		\item Define $W_{\varphi} = \{\Gamma_{D}\mid \exists D\text{ -- }\varphi\text{ -- description }: D\subseteq\Gamma_{D}\}$ and $|\psi|_{\varphi} = |\psi|_{\textup{L}}\cap W_{\varphi}$, for any $\psi\in Fm$. Clearly, $W_{\varphi}$ is finite.
	\end{enumerate}
	
	After that, Lewis proves the following lemma:
	
	\begin{lemma}\label{Lemma4}
		Let $\gamma$ be a truth-functional combination of subformulas of $\varphi$. Then $\gamma\in\Delta\text{ iff }\textup{L}\vdash\gamma$ for all $\Delta\in W_{\varphi}$.
	\end{lemma}
	
	Every element $\Delta$ in $W_{\varphi}$ has a unique $\varphi$ -- description as a subset. Let $\lambda_{D}$ be the first (in the enumeration of all formulas) conjunction of all the formulas in $D$. Lewis calls $\lambda_{D}$ a lable of $D$. For any subset $A\subseteq W_{\varphi}$ define lable of $A$ as follows: $\lambda_{A} = \bigvee\limits_{D\in A} \lambda_{D}$. Let also $\lambda_{\varnothing} = \bot$. Lewis defines $\bm{\Box}_{\varphi}$ in the following way:
	$$\bm{\Box} A = |\Box\lambda_{A}|_{\varphi},\text{ for }A\subseteq W_{\varphi}.$$ 
	
	Thus, we have model $M_{\varphi}=(W_{\varphi}, \bm{\Box}_{\varphi}, V_{\varphi})$, where $V_{\varphi}(p)=|p|_{\varphi}$.
	
	After this definition, Lewis tries to prove the following lemma, using labels. 
	
	\begin{lemma}\label{Lemma7}
		If $\psi$ is a Boolean combination of formulas from $Sub(\varphi)$, and $|\psi|_{\mathcal{M}_{\varphi}}=$ $|\psi|_{\varphi}$, then $|\Box \psi|_{\mathcal{M}_{\varphi}}=|\square \psi|_{\varphi}$.
	\end{lemma}
	
	Now we present the proof up to the part, where Lewis has an error.
	
	\begin{proof}
		Since $|\psi|_{\mathcal{M}_{\varphi}}=$ $|\psi|_{\varphi}$, we have that $$|\Box\psi|_{\mathcal{M}_{\varphi}} = \bm{\Box}|\psi|_{\mathcal{M}_{\varphi}}=\bm{\Box}|\psi|_{\varphi}=|\Box\lambda_{|\psi|_{\varphi}}|_{\varphi}.$$
		
		Now we have to show that $|\Box\lambda_{|\psi|_{\varphi}}|_{\varphi}=|\Box\psi|_{\varphi}$. By definition of labels we have $|\psi|_{\varphi}=|\lambda_{|\psi|_{\varphi}}|_{\varphi}$. Thus, $\psi\leftrightarrow\lambda_{|\psi|_{\varphi}}\in \Gamma_{D}$ for every $\Gamma_{D}$. By Lemma \ref{Lemma4} we have that $\textup{L}\vdash \psi\leftrightarrow\lambda_{|\psi|_{\varphi}}$. By the rule $\textup{(RE)}$ we have $\textup{L}\vdash \Box\psi\leftrightarrow\Box\lambda_{|\psi|_{\varphi}}$. And then Lewis applies Lemma \ref{Lemma4} to this statement to obtain $\Box\psi\leftrightarrow\Box\lambda_{|\psi|_{\varphi}}\in \Gamma_{D}$ for every $\Gamma_{D}$. However, this is not allowed, because Lemma \ref{Lemma4} can be applied to $\Box\psi\leftrightarrow\Box\lambda_{|\psi|_{\varphi}}$ only in case when both $\Box\psi$ and $\Box\lambda_{|\psi|_{\varphi}}$ are subformulas of $\varphi$. However, it may not be the case.
	\end{proof}
	
	Lemma \ref{Lemma7} is the key lemma for all the proof. However, it does not seem possible to fix it, because the error is forced by the construction of labels themselves.

\end{document}